\documentclass{amsart}                     

\usepackage{graphicx}
\usepackage{amsmath}
\usepackage{amssymb}
\usepackage{amsopn}
\usepackage{epsfig}
\usepackage{amsfonts}
\usepackage{latexsym}
\usepackage[english]{babel}

\newcommand{\expec}[1]{\ensuremath{{\rm E}\left[#1\right]}}
\newcommand{\var}[1]{\ensuremath{{\rm Var}\left( #1 \right)}}

\newcommand{\me}{\mathrm{e}}
\newcommand{\md}{\mathrm{d}}

\begin{document}

\title[Reactive solute transport.
From particle dynamics to the PDF's]{Understanding the non-Gaussian nature of reactive solute transport.
From particle dynamics to the partial differential
equations}


\author[G.~Uffink \and A.~Elfeki \and M.~Dekking \and J.~Bruining \and C.~Kraaikamp]{Gerard Uffink \and Amro Elfeki \and Michel~Dekking \and Johannes Bruining \and Cor~Kraaikamp}

\address{G.~Uffink and  J.~Bruining,
              Delft University of Technology,
              Department of Civil Engineering and Geosciences \\
              Stevinweg 1, 2628CN Delft, The Netherlands               
              }
\email{G.J.M.Uffink@tudelft.nl}                           
\address{A.~Elfeki,
              Dept.~of Water Resources,
Faculty of Metereorology, Environment and Arid Land Agriculture\\
King Abdul Aziz University\\ Jeddah, Saudi Arabia}
 \address{F.M.~Dekking and C.~Kraaikamp\\
              Delft University of Technology,
              EEMCS, Delft Institute of Applied Mathematics  \\
              Mekelweg 4, 2628CD Delft, The Netherlands}

\maketitle

\begin{abstract}
In the present study we examine non-Gaussian spreading of solutes
subject to advection, dispersion and kinetic sorption
(adsorption/desorption).  We start considering the behavior of a
single particle and apply a random walk to describe
advection/dispersion plus a Markov chain to describe kinetic
sorption. We show in a rigorous way that this model leads to a set
of differential equations. For this combination of stochastic
processes such a derivation is new. Then, to illustrate the
mechanism that leads to non-Gaussian spreading we analyze this set
of equations at first leaving out the Gaussian dispersion term
(microdispersion). The set of equations now transforms to the
telegrapher's equation. Characteristic for this system is a
longitudinal spreading, that becomes Gaussian only in the
long-time limit. We refer to this as kinetics induced spreading.
When the microdispersion process is included back again, the
characteristics of the telegraph equations are still present. Now
two spreading phenomena are active, the Gaussian microdispersive
spreading plus the kinetics induced non-Gaussian spreading. In the
long run the latter becomes Gaussian as well. Another non-Gaussian
feature shows itself in the 2D situation. Here, the lateral spread
and the longitudinal displacement are no longer independent, as
should be the case for a 2D Gaussian spreading process. In a
displacing plume this interdependence is displayed as a `tailing'
effect. We also analyze marginal and conditional moments, which
confirm this result. With respect to effective properties
(velocity and dispersion) we conclude that effective parameters
can  be defined properly only for large times (asymptotic times).
In the two-dimensional case it appears that the transverse
spreading depends on the longitudinal coordinate. This results in
`cigar-shaped' contours.

\end{abstract}

\section{Introduction}
\label{intro}It is well known that in the field a contaminant
plume spreads at a higher rate than as predicted by theory and lab
experiments. In addition, one observes that the spreading pattern
often deviates from the Gaussian pattern, especially at early
times (tailing). These phenomena have been analyzed intensively,
theoretically, numerically and by field experiments
(\cite{Gelhar}, \cite{Neuman}, \cite{Dagan}, \cite{Dieulin2},
\cite {Dieulin1}), \cite{Scheidegger}, \cite{Biggar}, \cite
{Maloszewski}). The majority of papers on this subject attribute
 these phenomena to the heterogeneity of the medium. In our paper we
show that in a homogeneous aquifer a similar behavior occurs when
contaminants are subject to a relatively slow kinetic
adsorption/desorption reaction.

We consider a homogeneous medium and follow a (solute) particle
during its motion through the pore system, while simultaneously it
is subject to adsorption/desorption. Although the porous medium
itself is simply homogeneous and non-stochastic, the particle's
behavior is still chaotic and, therefore, our particle model will
be stochastic. Note that our approach differs from that of papers
on particle tracking in the sense that we do not start with the
advection-dispersion-reaction equations and, accordingly, perform
a particle tracking algorithm to solve or simulate the equations.
Instead, we start with a stochastic description for the movement
of a single particle and work our way towards the differential
equations. Particle tracking papers assume an analogy between
particle tracking and the differential equations either by just
stating the validity of this analogy or by referring to previous
papers on particle tracking. To our knowledge a rigorous
derivation of the differential equations for advective-dispersive
transport with kinetic sorption starting from a stochastic model
for a single particle does not exist in any of those papers.

We also discuss particle models used in various other fields, such
as the velocity jump model (chemotaxis or movement of bacteria,
Hillel). This velocity jump model (Giddings-Eyring model) leads to
a telegraph equation. It differs from our model by the fact that
there is no diffusion or dispersion phenomenon. We will
demonstrate that, when the velocity jump is extended with
dispersion/diffusion, the telegraph character stills remains. This
explains the non-Gaussian behavior in the pre-asymptotic stage.

The non-Gaussian features can also be illustrated by (spatial)
moments, especially by the third (skewness) and fourth (kurtosis)
centralized moments. For the 1D situation first and second moments
were studied by previous authors e.g. \cite{ValJH},
\cite{Chrysikopoulos} and \cite{Souadnia}. These papers mainly
focused on the large time (asymptotic) results. Pre-asymptotic
expressions for first and second moments were first derived by
\cite{Michalak}. Unfortunately, in the formulation of the
solutions Michalak and Kitanidis suggest that for an arbitrary
initial distribution of the phases the mean and variance in each
separate phase can be obtained by a linear combination. However,
this is only true for the non-centralized moments, but not for the
variance. Moreover, one of the 4 solutions is incorrect (see
section \ref{1Dmoment}, also \cite{DekkingC}). In our paper we
shall consider the moments for the 2D situation and examine in
detail the conditional moments. These moments manifest more
clearly the non-Gaussian behavior than the 1D or marginal moments.
It shows that non-Gaussian behavior can exist, even when the
marginal moments suggest that the behavior is Gaussian.

In our `particle view' the behavior of a particle is described by
a random walk and a two state Markov process. A similar model was
used by \cite{ValJH}, \cite{Valocchi} and \cite{VanKoot}.
Equivalent processes occur in other fields, such as chromatography
(\cite{GidEyr}, \cite{Gid}, \cite{Keller}) or in statistical
physics (e.g. \cite{Furth}, \cite{Goldstein}, \cite{Kac},
\cite{Weiss} and \cite{Masoliver}). The amount of literature
related to this topic is extensive and spreads over many different
fields. Here, we look at this topic from the perspective of solute
transport in groundwater.

\section{The particle view}\label{sec:DT-CS}
\subsection{A stochastic particle model}\label{intro}
We model the movement  of a single particle subject to an
advection/dis\-persion/sorp\-tion process over a time interval
$[0,t]$. For simplicity we discuss the one-dimensional case. We
discretize time by choosing some $n$, and by dividing $[0,t]$ into
$n$ intervals of length  $\Delta t=t/n$. We observe the state of
the particle at the time points $0,\Delta t,2 \Delta t,
\dots,n\Delta t$. For this state there are two possibilities:
`free' or `adsorbed,' which we code by the letters $f$ and $a$.
The particle can only move when it is `free,' and its displacement
has two components: dispersion and advection.

Let $X_k$ be the displacement due to the dispersion of the
particle the $k$th time that it is `free.' We model the  $X_k$ as
independent random variables with mean and variance
\begin{equation*}
\expec{X_k}=0 \quad {\rm and } \quad   \var{X_k}=2D\Delta t.
\end{equation*}
The displacement due to advection is given by $v\Delta t$, where
$v$ is the (deterministic) advection velocity.

Let $K_n$ be the number of intervals $[k\Delta t,(k+1)\Delta t)$
during $[0,t]$ that the  particle was `free.' (Here the open
bracket at the right end indicates that the point $(k+1)\Delta t$
is not included.) In other words, $K_n\Delta t$ is the free
residence time of the particle in $[0,t]$. Let $S(t)$ be the
position of the particle at time $t$.  Combining the two types of
displacement we obtain
\begin{equation*}
S(t)=\sum_{k=1}^{K_n} (X_k +v\Delta t).
\end{equation*}
The distribution of $K_n$ is determined by the kinetics, i.e., by
the switching between the `free' and the `adsorbed' state. This is
naturally described by a two state Markov chain. The
state-transitions of this chain after a certain time step $\Delta
t$ are given by a transition probability matrix $(p_{ij})$:
\begin{equation*}
\left[ {\begin{array}{*{20}c}
   {p_{f\!f} } & {p_{f\!a} }  \\
   {p_{a\!f} } & {p_{a\!a} }  \\
\end{array}}  \right] =\left[ {\begin{array}{*{20}c}
   {1 - b} & {b}  \\
   {a} & {1-a}  \\
\end{array}} \right].
\label{kineticsmatrix}
\end{equation*}
This means that for instance the transition from `adsorbed' to
`free' has probability $a=p_{a\!f}$ to happen during the time
interval $[k\Delta t,(k+1)\Delta t)$ (note that actually we make
this change---if it takes place---at the end of the interval).

We will be mainly interested in the moments of $S(t)$. Below we
will compute the first and second moment, and in the next Section
we discuss the (centered) third and fourth moment. To compute
\expec{S(t)} we use the well known formula (see e.g.
\cite{Ross-Prob-Mod}) for a random sum of $K_n$ independent and
identically distributed random variables $Y_k$ (also independent
of $K_n$):
\begin{equation*}
\expec{\sum_{k=1}^{K_n} Y_k}=\expec{K_n}  \expec{ Y_1}.
\end{equation*}
Here the mean of $K_n$ equals:
$
\expec{K_n}=\frac{a}{a+b}\, n.
$
This expression can be obtained from~\cite{Viveros},
or~\cite{DekkingC}. Substituting we find (with $n\Delta t =t$)
\begin{equation*}
\begin{split}
 \expec{S(t)}&=\expec{ K_n} \left(\expec{
X_1+v\Delta t}\right)= \frac{a}{a+b}\, v t.
\end{split}
\label{Retardation}
\end{equation*}
Here,  $(a+b)/a$ is  the retardation factor $R$. To see this, note
that the probability vector $\big(b/(a+b)\quad a/(a+b)\big)$ is
the stationary distribution of the Markov chain, and so $a/(a+b)$
is the fraction of time the particle is free. Thus the effective
velocity $v^*=v/R$. We have implicitly required that the particle
at time $0$ is given the state `adsorbed' or `free' according to
this distribution, for other initial distributions there will be a
correction term in the formula for \expec{S(t)}, which tends to
$0$ as $t$ tends to infinity (cf. \cite{DekkingC}).

To compute \var{S(t)} we use  the well known formula (see e.g.
\cite{Ross-Prob-Mod}) for the {\it{variance}} of a random sum of
$K_n$ i.i.d. random variables $Y_k$ (also independent of $K_n$):
\begin{equation}\label{var-random-sum}
\var{ \sum_{k=1}^{K_n} Y_k} = \expec{ K_n}
\var{Y_1}+\var{K_n}(\expec{ Y_1})^2.
\end{equation}
This yields with $Y_k=X_k +v\Delta t$ and $n\Delta t =t$:
\begin{equation}\label{vardiscrete}
\begin{split}
\var{S(t)}  =  & \,\expec{ K_n} \var{X_k+v\Delta t}
               +\var{K_n}(\expec{X_k+v\Delta t})^2\\
            =  & \frac{a}{a+b}\,2D\, t+ \var{K_n} v^2(\Delta t)^2,
\end{split}
\end{equation}
where (as can be deduced from \cite{Viveros} or \cite{DekkingC})
\begin{equation}
\label{varK_n}
\var{K_n}=\frac{ab(2-a-b)}{(a+b)^3}\,n
-\frac{2ab(1-a-b)}{(a+b)^4}[1\!-\!(1\!-\!a\!-b)^n].
\end{equation}
The equations (\ref{vardiscrete}) and (\ref{varK_n}) thus tell us
that the variance of the displacement of the particle grows more
or less linearly in time with (asymptotic) slope
\begin{equation*}
\dfrac{a}{a+b}2D+ \dfrac{ab(2-a-b)}{(a+b)^3} v^2\Delta t.
\end{equation*}

\subsection{Skewness and kurtosis}
\label{MandV} To obtain the skewness of $S(t)$ we must use  the
not so well known formula for the {\it{third central moment}} of a
random sum of $K_n$ i.i.d. random variables $Y_k$ (also
independent of $K_n$):
\begin{equation*}
\begin{split}
&\expec{\big(S(t)-\expec{S(t)}\big)^3} =\expec{\Big(\sum_{k=1}^{K_n} Y_k-\expec{ \sum_{k=1}^{K_n} Y_k }\Big)^3}\\
& =\expec{K_n} \expec{(Y_1-\expec{Y_1})^3} + 3\expec{Y_1}\var{Y_1}\var{K_n}\\
&\qquad + (\expec{Y_1})^3\,\expec{(K_n-\expec{ K_n})^3}.
\end{split}
\end{equation*}
To actually derive a formula for the skewness of the displacement
of the particle from this equation will lead to very heavy
computations (and the situation for the kurtosis is even worse).
However, without doing any computations we can already tell that
as $t\rightarrow\infty$ the skewness must tend to zero, and the
kurtosis to 3: this is because the distribution of  the
displacement of the particle will tend to a Gaussian by the
Central Limit Theorem for \emph{random} sums of independent
identically distributed random variables. In our case this follows
since $K_n/n$ tends in the mean, and hence in probability to
$a/(a+b)$, see \cite{Feller}, page 258.

\subsection{Decreasing the time steps}
The discrete time steps are somewhat unnatural.
 We would like to let $\Delta t$ tend to 0.
But then we have to realize that $a$ and $b$ are functions of
$\Delta t$. Since the probability that the particle changes its
state is proportional to the time $\Delta t$ it is observed (if
$\Delta t$ is not too large), we should put
\begin{equation*}
a=\mu \Delta t,\qquad b=\lambda \Delta t, \label{amu}
\end{equation*}
where $\mu$ and $\lambda$ are now the \emph{rates} at which the
particle switches from `adsorbed' to `free', and from  `free' to
`adsorbed'.
 Substituting this in
Eqs~(\ref{vardiscrete}) and (\ref{varK_n}), we obtain
\begin{equation*}
\begin{split}
 \var{S(t)}   &= \frac{\mu}{\lambda+\mu}2Dt+
         \frac{\lambda\mu(2-(\lambda+\mu) \Delta t)}{(\lambda+\mu)^3}   v^2t\\
      &\qquad   -\frac{2\lambda\mu(1-(\lambda+\mu) \Delta t)}{(\lambda+\mu)^4}
        \left[1-{\left(1-(\lambda+\mu)\frac{t}{n}\right)\!}^n\right] v^2.
\end{split}
\end{equation*}
Letting $\Delta t\rightarrow 0$, and hence $n\rightarrow\infty$ we
obtain
\begin{equation*}
\begin{split}
\var{S(t)} = \frac{\mu}{\lambda+\mu}2Dt+
         \frac{2\lambda\mu}{(\lambda+\mu)^3} v^2t
     -\frac{2\lambda\mu}{(\lambda+\mu)^4}\left[1-\me^{-(\lambda+\mu)t}\right] v^2.
\end{split}\label{variance-total}
\end{equation*}

Thus we recuperate a (more general and more detailed) version of
the main result of \cite{gut}, \emph{and} there is  a match with
the expression that comes from the moment analysis based on the
differential equations (as can be derived by correcting the
results in \cite{Michalak}).
\subsection{The state of the particle at time $t$}
\label{MandV-totalplume2} Let $S_f(t)$ be the position of the
particle at time $t$ \emph{given} that it is free at time $t$.\\
To find the distribution of $S_f(t)$, we need the distribution of
$K_n^{(f)}$, the number of intervals $[k\Delta t,(k+1)\Delta t)$
during $[0,t]$ that the  particle was free, given that it is free
at time $t=n\Delta t$. We find now (where \expec{K_n^{(f)}} can be
deduced from \cite{Viveros}) that:
\begin{equation*}\label{mean-free}
\begin{split}
\expec{S_f(t)}&=\expec{K_n^{(f)}}v\Delta t
= \left[ \frac{a}{a+b}\,n + \frac {b\,(1-(1 - a -b)^n) }{(a + b)^{2}}\right] v\Delta  t\\
    &=\frac{a}{a+b}\, v t + \frac {b\,(1-(1 - a -b)^{t/\Delta t}) }{(a + b)^2}\, v \Delta t.
\end{split}
\end{equation*}
Substituting $a=\mu \Delta t,\, b=\lambda \Delta t,$ and letting
$\Delta t\rightarrow 0$ we obtain
\begin{equation*}\label{mean-free}
\expec{S_f(t)}= \frac{\mu}{\lambda\!+\!\mu}\,vt +\frac{\lambda
v}{(\lambda+\mu)^2}\left(1\!-\!\me^{-(\lambda\!+\!\mu)t}\right).
\end{equation*}
From \cite{DekkingC} we have that $\var{K_n^{(f)}}$ equals
\begin{equation*}
\begin{split}
& \left[{\displaystyle \frac {a\,b\,(2 - a -b) }{(a + b)^{3}}}
  + {\displaystyle \frac{2\,b\,(a - b)\,(1 - a- b)^{n}}{(a + b)^{3}}}\right]\,n\\
 &  \qquad + \left[ \frac{b\,(3\,a -b)}{(a + b)^3} - \frac{4\,a\,b}{(a + b)^4} \right][1-(1 - a
 -b)^n]+ \frac {b^2}{(a +b)^4}\,[1-(1 - a - b)^{2\,n}].
\end{split}
\end{equation*}
Using Equation (\ref{var-random-sum}) we derive from this
\begin{equation*}
\begin{split}
\var{S_f(t)} & =   \expec{ K_n^{(f)}} \var{X_k+v\Delta t}
             +\var{K_n^{(f)}}(\expec{X_k+v\Delta t})^2\\
          & =\frac{a}{a+b}\,2Dt +\var{K_n^{(f)}} v^2(\Delta t)^2
          + \frac {b\,(1-(1 - a -b))^{t/\Delta t}}{(a + b)^2}2D\Delta t.
\end{split}
\end{equation*}
Substituting $a=\mu \Delta t, b=\lambda \Delta t$, and letting
$\Delta t \rightarrow 0$ we obtain
\begin{equation*}\label{variance-free}
\begin{split}
\var{S_f(t)} & =  \frac{\mu}{\lambda+\mu}\,2Dt
+\frac{\lambda}{(\lambda+\mu)^2}\,\left(1-\me^{-(\lambda+\mu)t}\right)2D\\
&\quad +  \left[ \frac{2\lambda\mu}{(\lambda+\mu)^3}+ \frac{2\lambda(\mu-\lambda)}{(\lambda+\mu)^3}\,\me^{-(\lambda+\mu)t}\right]\,v^2t\\
&\qquad -\frac{4\lambda\mu}{(\lambda+\mu)^4}\,\left[1-\me^{-(\lambda+\mu)t}\right]v^2+
\frac{\lambda^2}{(\lambda+\mu)^4}\,\left[1-\me^{-2(\lambda+\mu)t}\right]v^2.
\end{split}
\end{equation*}
It can be shown that this matches with the expressions in
\cite{Michalak}, when these are corrected as in \cite{DekkingC}.
Similar computations can be made for the displacement of the
particle \emph{given} that it is absorbed at time $t$.
\subsection{Derivation of the differential equations}
\label{sec:DT-DS} We will now show how the fundamental
differential equations~(\ref{Eq1dfree}),~(\ref{Eq1dads}) can be
obtained from a diffusion limit of the single particle model. Our
approach is similar to the one followed for transport in fluidized
beds in \cite{dehling}. In order to obtain this diffusion limit we
also discretize space in locations
\begin{equation*}
 i\Delta x, \quad i=\dots,-1,0,1,\dots.
\end{equation*}
Here we let $\Delta x$ depend on $\Delta t$ in the classical way,
which is motivated by the fact that typically at time $t$ the
spatial fluctuations are of order $\sqrt{t}$:
\begin{equation}
\label{scaling} \Delta x=c\sqrt{\Delta t},
\end{equation}
 where $c>0$ will be
chosen later. The particle moves according to a Markov chain
$(Z_n)$, which is a birth-death process (birth=one step to the
right, death=one step to the left), with the additional
possibility that the particle may become adsorbed and free again.
The state space is therefore a   product
\begin{equation*}
S=\{\dots,-1,0,1,\dots\}\times\{a,f\},
\end{equation*}
where e.g.~$Z_n=(i,a)$ means that at time $n\Delta t$ the particle
is at $ i\Delta x$, and is absorbed.
\begin{figure}[h!]\centering
\includegraphics[height=3cm]{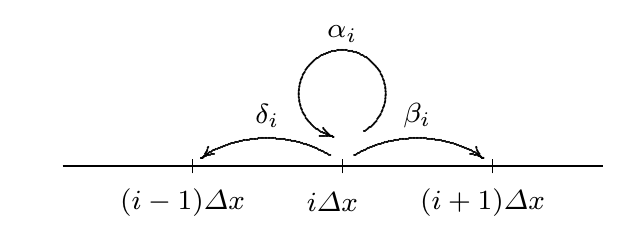}
\caption{The three possible movements of the particle when it is
free.} \label{fig:free-move}
\end{figure}
Let $\beta_i, \delta_i $ and $\alpha_i$ be the probabilities that
the  particle (in the free state) moves  from  $ i\Delta  x$  to $
(i+1)\Delta x$, from $i \Delta x$ to  $ (i-1)\Delta x$, or stays
at $ i\Delta  x$. Note that in this model the dispersion may
depend on the location, i.e., we could more generally consider
(deterministic) inhomogeneous media. To make this more realistic
we could also let the adsorption and desorption rates  depend on
$i$. As we consider here only the homogeneous case, we take
$\beta_i=\beta, \delta_i=\delta $ and $\alpha_i=\alpha$ for each
$i$. To fit the new model  with the  model in Section \ref{intro}
we    choose values for $\beta, \delta $ and $\alpha$ in the
following fashion:
\begin{equation*}
\begin{split}
\beta & = D\frac{\Delta t}{(\Delta x)^2}\,+
    \,\frac{v^2}{2}\Big(\frac{\Delta t}{\Delta  x}\Big)^2+\,\frac{v}{2}\frac{\Delta t}{\Delta x}\\
\delta & = D\frac{\Delta t}{(\Delta x)^2}+
    \,\frac{v^2}{2}\Big(\frac{\Delta t}{\Delta  x}\Big)^2-\,\frac{v}{2}\frac{\Delta t}{\Delta x}\\
\alpha & = 1- \beta-\delta \, =\, 1- 2D\frac{\Delta t}{(\Delta
x)^2}\,-\,
    v^2\Big(\frac{\Delta t}{\Delta  x}\Big)^2.
\end{split}
\end{equation*}
By fitting we mean that the mean and the variance of the
displacement in a time interval of length $\Delta t$ of the
particle (in the free state)  are the same in the two models.
Indeed, the mean of this displacement equals $\beta\Delta x
-\delta \Delta x = v\Delta t$, and the variance equals
\begin{equation*}
\beta(\Delta x)^2 +\delta (\Delta x)^2 -(v\Delta t)^2= 2D \Delta
t.
\end{equation*}
In terms of $\Delta t$ only, using (\ref{scaling}), the
displacement probabilities are
\begin{equation*}
\begin{split}
\beta  &=  \frac{D}{c^2}+\frac{v^2\Delta
t}{2c^2}+\frac{v\sqrt{\Delta t}}{2c} \\
\delta  &= \frac{D}{c^2}+\frac{v^2\Delta
t}{2c^2}-\frac{v\sqrt{\Delta t}}{2c} \\
\alpha  &= 1-\frac{2D}{c^2}-\frac{v^2\Delta t}{c^2}.
\end{split}
\end{equation*}
From this we see that these are indeed probabilities for $\Delta
t$ small enough, provided we choose $c> \sqrt{2D}$. The possible
transitions of the chain are

\begin{figure}[h!]\centering
\includegraphics[height=4cm]{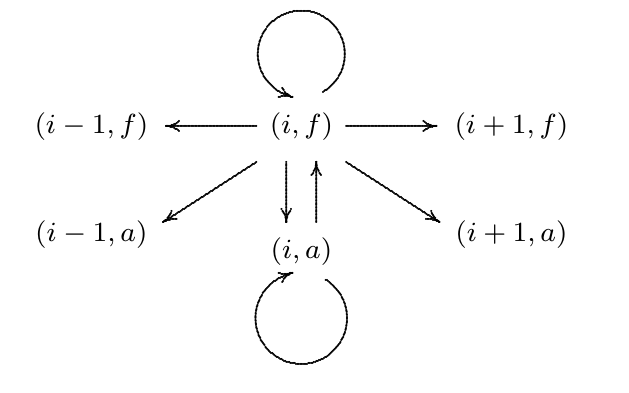}
\end{figure}

The corresponding transition probabilities  are:
$$
\begin{array}{lll}
p_{(i,f),(i+1,f)} = \beta(1-\lambda \Delta t) \quad &
p_{(i,f),(i+1,a)} = \beta \lambda \Delta t \quad &
p_{(i,f),(i-1,f)}  =  \delta (1-\lambda \Delta t)\\
p_{(i,f),(i-1,a)} = \delta \lambda \Delta t &
p_{(i,a),(i,a)} =  1-\mu \Delta t & p_{(i,a),(i,f)} = \mu
\Delta t\\
p_{(i,f),(i,f)} =  \alpha (1-\lambda \Delta t) &
p_{(i,f),(i,a)} =  \alpha \lambda \Delta t &
\end{array}
$$
Let $p_f(i,n)$, respectively $p_a(i,n)$, be the probability that
at time $t=n \Delta t$ the particle is free, respectively absorbed
and at position $i\Delta x$.

\noindent The master equations for the particle are
\begin{equation} \label{mobile master}
\begin{split}
p_f(i,n+1) & = p_f(i-1,n)\beta(1-\lambda\Delta t)
+ p_f(i+1,n)\delta(1-\lambda\Delta t)\\
    & \qquad + p_f(i,n)\alpha(1-\lambda\Delta t)  + p_a(i,n)\mu\Delta t,
\end{split}
\end{equation}
and
\begin{equation}
\begin{split}  \label{absorbed master}
p_a(i,n+1) &    =p_f(i-1,n)\beta \lambda\Delta t+p_f(i+1,n)\delta\lambda\Delta t\\
    & \qquad +p_f(i,n)\alpha \lambda\Delta t  + p_a(i,n)(1-\mu\Delta t).
\end{split}
\end{equation}
With some abuse in notation regarding the functions $p_a$ and
$p_f$ we will denote the limiting probabilities of $p_a(i,n)$ and
$p_f(i,n)$ as $i\Delta x\rightarrow x$, and $t=n\Delta t$ by
$p_a(x,t)$  and $p_f(x,t)$. We will obtain partial differential
equations for these limiting probabilities   when we let $\Delta
t\rightarrow 0$,  and $i=i(\Delta t)\rightarrow \infty$, in such a
way that  $i\Delta x\rightarrow x$. (The obvious way to achieve
this is to take $i(\Delta t)$ equal to the integer closest  to
$x/(c\sqrt{\Delta t})$.) Rearranging (\ref{mobile master}) we
obtain
\begin{equation}
\begin{split}\label{mobile-rec1}
p_f(i,n+1) &
  = p_f(i-1,n)\beta+p_f(i+1,n)\delta
+p_f(i,n)\alpha
   -\lambda\Delta t [p_f(i-1,n)\beta\\
&\qquad +p_f(i+1,n)\delta
    +p_f(i,n)\alpha]
+ p_a(i,n)\mu\Delta t.
\end{split}
\end{equation}
In the first three terms of the right hand side we substitute the
values for $\beta, \delta $ and   $\alpha$:
\begin{equation}
\begin{split}
\label{mobile-rec2}
p_f(i-1,n)\beta & +p_f(i+1,n)\delta  +p_f(i,n)\alpha = \\
   & = \Big[D\frac{\Delta t}{(\Delta x)^2}\,+
        \,\frac{v^2}{2}\Big(\frac{\Delta t}{\Delta  x}\Big)^2\Big]
            \Big[p_f(i-1,n)-2p_f(i,n)
            +p_f(i+1,n)\Big]\\
         &\qquad + \Big[\frac{v}{2}\frac{\Delta t}{\Delta x}\Big]
                    \Big[p_f(i-1,n)+p_f(i+1,n)\Big]+ p_f(i,n).
\end{split}
\end{equation}
Substituting Equation (\ref{mobile-rec2}) into (\ref{mobile-rec1})
and dividing by $\Delta t$, we obtain an equation for the
difference quotient $(p_f(i+1,n)-p_f(i,n))/\Delta t$. Then,
letting $\Delta t\rightarrow 0$, we obtain
\begin{equation}\label{Eq1dfree}
\frac{{\partial p_f(x,t) }}{{\partial t}} =  D\frac{{\partial ^2
p_f(x,t) }}{{\partial x^2 }} - v\frac{{\partial p_f(x,t)
}}{{\partial x}}
 - \lambda p_f(x,t)  + \mu p_a(x,t).
\end{equation}
For the adsorbed phase a similar equation can be derived. Here,
since the adsorbed particle is not subject to advection or
dispersion, the first and second term at the right hand side are
absent.
\begin{equation}
\frac{{\partial p_a(x,t) }}{{\partial t}} = - \mu p_a(x,t) +
\lambda p_f(x,t). \label{Eq1dads}
\end{equation}
\subsection{From particle to plume}
It might seem surprising that we study the behavior of a
contaminant plume from the stochastic analysis of a single
particle. Here we illustrate how these two approaches are
connected. We model the contaminant plume by a collection of $N$
particles. These particles move independently according to the
same law as the single particle considered in the previous
sections.
 Let $S_i(t)$ be
the position of the $i$th particle at time $t$. We are interested
in the centroid $Z(t)$ of the plume at time $t$. This is given  by
\begin{equation*}
Z(t)=\frac1{N}\sum_{i=1}^N S_i(t).
\end{equation*}
We are also interested in the spreading of the plume around its
centroid. This we measure by the (empirical) variance $V(t)$ of
the particles given by
\begin{equation}\label{plumevariance}
 V(t)=\frac1{N}\sum_{i=1}^N \left(S_i(t)-Z(t)\right)^2.
\end{equation}

The random variable $Z(t)$ is an average of independent
identically distributed random variables with finite expectation.
Therefore by the strong law of large numbers for $t$ fixed, and
$N$ large
\begin{equation}\label{centroidapprox}
 Z(t)\approx \expec{S(t)}.
\end{equation}

We now turn to the spread of the plume. Here the situation is more
complicated because the terms in the sum are no longer independent
random variables, and also depend on $N$. It is well known that a
rewriting of Equation (\ref{plumevariance}) yields
\begin{equation*}
\begin{split}
V(t)&=\frac1{N}\sum_{i=1}^N \left(S_i(t)-Z(t)\right)^2
= \frac1{N}\sum_{i=1}^N
\left(S_i(t)-\expec{S(t)}\right)^2\!-\!(Z(t)\!-\!\expec{S(t)})^2.
\end{split}
\end{equation*}
We have already argued (see (\ref{centroidapprox})) that the last
term is approximately 0, and another application of the law of
large numbers to the first term yields that
\begin{equation*}
V(t)\approx \expec{(S(t)-\expec{S(t)})^2} =\var{S(t)}.
\end{equation*}
\section{Giddings-Eyring Model}
\label{Analytical Solutions} A simpler system, different from but
still related to the system described by Eqs. (\ref{Eq1dfree}) and
(\ref{Eq1dads}), is given by the following equations:
\begin{eqnarray}
\frac{{\partial N_f(x) }}{{\partial t}} + v\frac{{\partial N_f(x)
}}{{\partial x}} &=& \mu  N_a(x)  - \lambda  N_f(x)\nonumber \\
 \frac{{\partial N_a(x)
}}{{\partial t}} &=& \lambda  N_f(x)  - \mu  N_a(x).\nonumber
\end{eqnarray}
The difference is that the dispersion process is absent. This
system is known in the literature under various names. In
probability theory and statistical physics it is known as a
persistent or correlated random walk and has been studied e.g.\ by
\cite{Furth}, \cite{Goldstein}, \cite{taylor}, \cite{Kac}, and
\cite{Weiss}. In the field of chromatography the system is
intensively studied as well (\cite{GidEyr}, \cite{Gid},
\cite{Keller}). Solutions for the probability density functions of
the particles are given by Giddings and Eyring. \cite {HillelT},
cite{Hillel} applies the equations to the movement of bacteria in
the direction of the gradient of food molecules (chemotaxis) and
uses the term ``velocity jump process''. Several interesting
observations can be made with respect to the spreading of the
particles. The particles are undergoing an `apparent' dispersion,
despite the fact that (hydrodynamic) dispersion is not included in
the model. This `kinetics-induced' dispersion develops in a
non-Gaussian way. In systems where both hydrodynamic dispersion
and kinetics-induced dispersion are present, the latter sometimes
can be more dominant, such that the non-Gaussian spreading is
observed also in systems with hydrodynamic dispersion. We shall
illustrate this later. First, we introduce a moving coordinate
system with velocity $v^*$. The new $x$-coordinate is:
\begin{equation*}
\bar{x} =x-v^*t,
\end{equation*}
where $v^* = v \mu /(\lambda + \mu)$. In the new coordinate system
free particles move to the right with velocity $v_f = v - v^* = v
\lambda /(\lambda +\mu )$. The adsorbed particles `move' with
velocity $v_a = - v^* = -v\times\mu /(\lambda +\mu )$. The minus
sign indicates that the movement is to the left (i.e. with respect
to the new coordinate system). The equations are now:
\begin{equation}
\begin{split}
\frac{{\partial N_f(\bar x) }}{{\partial t}} + v_f\frac{{\partial
N_f(\bar x) }}{{\partial \bar{x}}} &= \mu  N_a(\bar x)  - \lambda
N_f( \bar x)\nonumber \\ \frac{{\partial N_a(\bar x) }}{{\partial
t}} + v_a\frac{{\partial N_a(\bar x) }}{{\partial \bar{x}}} &=
\lambda  N_f(\bar x)  - \mu  N_a (\bar x).\nonumber
\end{split}
\end{equation}
The equations can be rewritten as a single differential equation
by considering first the sum and difference of the particle
distributions (Kac's trick), i.e.,
\begin{equation*}
u(x) = N_f(x)  + N_a(x); \quad w(x) = N_f(x)  - N_a(x).
\end{equation*}
After summation and substraction of the differential equations we
obtain for $u$ and $w$:
\begin{eqnarray*}
\frac{{\partial u}}{{\partial t}}+
    \left( {\frac{{\lambda-\mu}}{{\lambda+\mu}}} \right)
        \frac{v}{2}\frac{{\partial u}}{{\partial\bar x}}+
        \frac{v}{2}\frac{{\partial w}}{{\partial \bar x}} &=& 0 \\
\frac{{\partial w}}{{\partial t}}
    + \left( {\frac{{\lambda-\mu }}{{\lambda+\mu }}}\right)
        \frac{v}{2}\frac{{\partial w}}{{\partial \bar x}}+
        \frac{v}{2}\frac{{\partial
u}}{{\partial \bar x}}&=&- u(\lambda-\mu)     -w(\lambda+\mu).
\end{eqnarray*}
Now we differentiate the first equation to $t$, the second to
$\bar x$ and eliminate the derivatives of $w$:
\begin{equation}
\begin{split}
 {\frac{{\lambda  \mu  v^2 }}{{\left( {\lambda   + \mu  }
\right)^3 }}} \frac{{\partial ^2 u}}{{\partial \bar{x}^{2} }} -
{\frac{1}{{\lambda   + \mu  }}} \frac{{\partial ^2 u}}{{\partial
t^2 }} - {\frac{{v\left( {\lambda  - \mu } \right)}}{{\left(
{\lambda  + \mu  } \right)^2 }}} \frac{{\partial ^2 u}}{{\partial
t\partial \bar{x}}} = \frac{{\partial u}}{{\partial t}}
\label{telegraph1}
\end{split}
\end{equation}
This is a telegrapher's equation with an additional term due to
asymmetry ($ \mu  \ne \lambda$). Hillel discusses the symmetrical
($\lambda = \mu$) case only, but Weiss (\cite{Weiss}) also makes
some remarks on asymmetry. Also see Masoliver \cite {Masoliver}
and ~\cite{Chandra}. The telegraph equation may be interpreted
either as a diffusion equation with a perturbation term that
disappears at large times, or as a wave equation with a
perturbation term that disappears at early times. Thus, as time
proceeds, the system can be described by three different
equations: first, a wave equation for early times; secondly, a
telegraph equation for intermediate times and thirdly, an
advection-dispersion equation for large times.

\subsection{Large time behavior and the advection-dispersion equation} Hillel uses the following argument that
leads to a useful result for large times. For large times the
velocity $[LT^{-1}]$ and sorption rates $[T^{-1}]$ typically are
expressed in large time units and thus their values become large.
Then, the first term at the left hand side of (\ref{telegraph1})
dominates over the second and third term. Therefore, at large
times the equation approximately describes a dispersion process
with an equivalent dispersion coefficient $D^*$, purely induced by
the kinetics:
\begin{equation*}
D^*  = \frac{{\lambda \mu v^2 }}{{(\lambda  + \mu )^3 }}.
\label{effectivedis}
\end{equation*}
\subsection{Short time behavior and the wave equation} In a similar way it can be shown
that for small times the terms at
 the left hand remain, while the right side becomes small.
 The remaining expression is a wave equation that can be written
 as:
\begin{equation*}
\left[ {\frac{\partial }{{\partial t}} + \frac{{\mu v}}{{\lambda +
\mu }}\frac{\partial }{{\partial \bar x}}} \right]\left[
{\frac{\partial }{{\partial t}} - \frac{{\lambda v}}{{\lambda  +
\mu }}\frac{\partial }{{\partial \bar x}}} \right]u = 0
\end{equation*}
For an initial pulse at $\bar x = 0$ the solution represents two
pulses propagating along the characteristics. In $\bar x$-space,
\[
\begin{array}{l}
\bar x  - v_a t  = 0 ,\qquad  \bar x - v_f t  = 0 \\
 \end{array}
 \]
which in $x$-space correspond to:\quad $x = 0$, $x - v t = 0$.

These pulses can be identified as a stagnant pulse (the original
adsorbed particles) and a travelling pulse (original free
particles).

\subsection{Intermediate time behavior}

The wave equation and diffusion equation are approximations for
the process at short and large times respectively. The
intermediate time is described exactly by the full telegraphers's
equation. Therefore, examination of the solutions to this equation
will give information on the pre-asymptotic spreading behavior of
this process.

Solutions for the distribution of the particles have been derived
by \cite{GidEyr}, \cite{Gid}, and~\cite{Keller}. For detailed discussions of these
functions see \cite {Genucht}, \cite {Lassey2} and \cite
{VanKoot}.

Giddings and Eyring consider four types of densities  $h_{f\!f}
,h_{a\!f} ,h_{f\!a} ,h_{aa}$:
\begin{equation}\begin{split}
h_{f\!f} (\tau ,t)& = \me^{ - \lambda \tau  - \mu (t-\tau)} \sqrt
{\frac{{\lambda \mu \tau }}{{t - \tau }}} I_1 \left( \theta
\right)+\me^{-\lambda t}\delta(t-\tau)\\ h_{f\!a} (\tau ,t)& =
\lambda \me^{ - \lambda \tau - \mu (t-\tau) } I_0 \left( \theta
\right)\\ h_{a\!f} (\tau ,t)& =  \mu \me^{ - \lambda \tau - \mu
(t-\tau) } I_0 \left(
\theta \right)\\
 h_{a\!a} (\tau ,t)& =
\me^{ - \lambda \tau  - \mu (t-\tau)} \sqrt {\frac{{\lambda \mu
(t-\tau) }}{\tau }} I_1 \left( \theta \right)+\me^{-\mu
t}\delta(\tau),
\end{split}\label{Gidtau}\end{equation}

\noindent where  $ \tau = x/v$, $\theta = {2\sqrt {\lambda \mu
\tau (t - \tau )} } $  and $I_0 (\cdot)$ and $I_1 (\cdot)$ are
modified Bessel functions. Note that $ \tau = x/v$ is not simply a
convenient scaling of the $x$-coordinate, but $\tau$ also
represents the residence time in the free phase. The expressions
$h_{ij}$ represent the probability densities of the free residence
time for different phases and different initial states of the
particles. The first index indicates the initial state of the
particle and the second index indicates the state of the particles
the pdf is referring to. The distributions are zero for $\tau < 0$
and $\tau
> t$ ($x < 0$ and $x > v t$). The delta functions at $\tau = 0$
and  $\tau = 0$ (which in $x$-space corresponds to $x  =  0 $ and
$x = v t$) represent exponentially decreasing
 pulses and can also be identified as the fractions of particles
that, since $t = 0$, did not (yet) perform a change of state. 

\begin{figure}[h]
\begin{center}
\includegraphics[width =10cm]{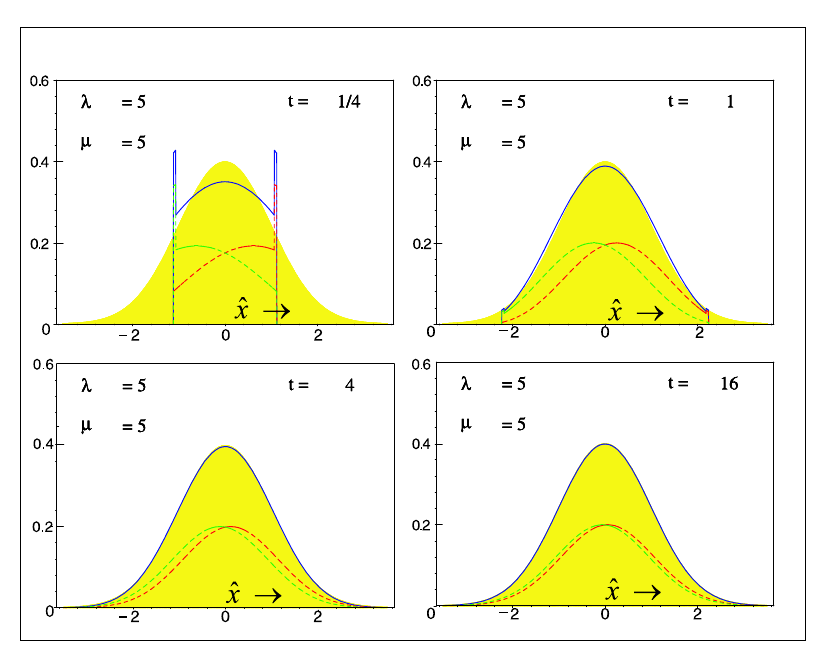}
\caption {Distributions of separate phases and total plume at 4
moments in time. Symmetric case ($\lambda = \mu$). Solid (blue)
line for total plume. Dashed lines for (green) adsorbed phase and
(red) free phase, yellow background for Gaussian
distribution.}\label{figureno1}
\end{center}
\end{figure}

\begin{figure}[t]
\begin{center}
\includegraphics[width =10cm]{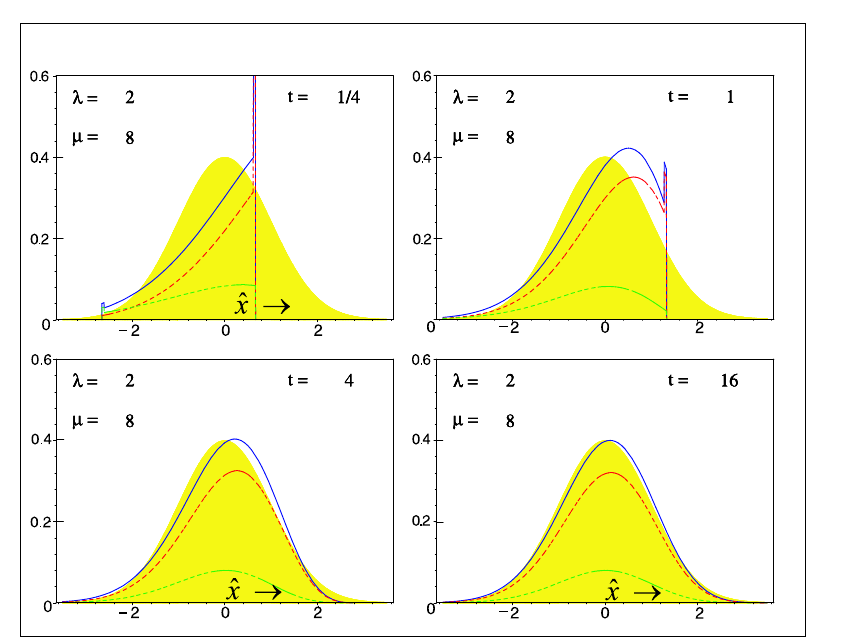}
\caption {Distributions of separate phases and total plume at 4
moments in time. Asymmetric case ($\lambda \neq \mu$). Solid
(blue) line for total plume. Dashed lines for (green) adsorbed
phase and (red) free phase, yellow background for Gaussian
distribution.}\label{figureno2}
\end{center}
\end{figure}

In
Figure \ref{figureno1} and \ref{figureno2} we present graphs with
the evolution in time of these distributions, using as initial
condition the (unit) pulse consisting of free and adsorbed
particles in equilibrium. Let $\pi_f$ and $\pi_a$ be the initial
amount of particles in each phase. Equilibrium exists for
$\frac{\pi_f}{\pi_a}= \frac{\mu}{\lambda}$, and, if the total
amount is unity we have:
\begin{equation*}
\pi_a = \frac{\mu}{\lambda+\mu}; \quad \pi_f =
\frac{\lambda}{\lambda+\mu}.\label{initiallyequilibrium}
\end{equation*}

Thus, for $t=0$:
\begin{equation*}
N_f (x,t) = \pi_f \delta (x);\quad N_a (x,t) = \pi_a \delta (x).
\end{equation*}
We denote the residence distributions due to this initial
condition as $h_f^{eq}$ and $h_a^{eq}$:
\begin{equation}
\begin{split}
h_f^{eq}(\tau,t) = \pi_f h_{f\!f}(\tau,t) +\pi_a h_{a\!f}(\tau,t)\\
h_a^{eq}(\tau,t) = \pi_f h_{f\!a}(\tau,t) + \pi_a
h_{a\!a}(\tau,t)\label{eq_definition},
\end{split}
\end{equation}
 while for the total amount of
particles we have, $h_{tot}^{eq} = h_f^{eq} + h_a^{eq}$.

In Figures \ref{figureno1} and \ref{figureno2} we compare
distributions (\ref{eq_definition}) and  $h_{tot}^{eq}$ with a
Gaussian distribution with mean $v^*t$ and variance $2D^*t$
representing the distribution of a solute with velocity $v^*$ and
dispersion $D^*$. All distributions are plotted versus a scaled
variable $\hat x$:
\begin{equation*}
\hat x = \frac{{x - v^* t}}{{\sqrt {2D^{^*} t} }}.
\end{equation*}
Several stages of the system are shown. At early times ($t = 1/4$,
Figure \ref{figureno1}, upper left graph) the pulses of free and
adsorbed mass are still distinguishable. The pulses move apart in
$\hat x$ space and are damped. The mass that `leaves' the pulses
gradually fills the space in between and builds up a distribution
that becomes Gaussian in the end (e.g., $t = 16$,
Figure~\ref{figureno1}, lower right graph). The distribution in
the interval between the pulses is absent in a pure wave system
and is typical for the telegraph equation. In the final stage the
pulses are completely damped and the distribution approaches the
normal distribution. Summarized, at early times the
`wave-character' dominates, at large times the
`diffusion-character' dominates, while at intermediate times the
system is adequately described by a telegraph equation (travelling
and dampened pulses + mixed zone in between). For $\lambda = \mu$
the Gaussian distribution is reached slightly faster than in the
asymmetric case ($\lambda \ne \mu$), as seen in Figure
\ref{figureno2}. For both cases the pulses disappear for $t$
exceeding both $3/\lambda$ and $3/\mu$. In the next section local
dispersion is included. We show that the two-dimensional
distribution still may deviate from two-dimensional Gaussian
functions even though the corresponding one-dimensional
distribution is close to Gaussian.

\section{Giddings-Eyring model including dispersion} We
extend the Gidding-Eyring model by including longitudinal and
transverse dispersion. This way we obtain the following
two-dimensional system.
\begin{equation}
\begin{split}
\frac{{\partial N_f }}{{\partial t}} \!-\! D_L \frac{{\partial ^2
N_f }}{{\partial x^2 }} \!-\! D_T \frac{{\partial ^2 N_f
}}{{\partial y^2 }} \!+\! v\frac{{\partial N_f
}}{{\partial x}} &= \mu  N_a \!-\! \lambda  N_f\\
 \frac{\partial N_a}{\partial t} &= \lambda  N_f  \!-\! \mu  N_a,\label{eq2D-2}
\end{split}
\end{equation}
where now $N_f$ and $N_f$ represent a two dimensional particle
density distribution $ \left[ {L^{ - 2} } \right] $ and $D_L$ and
$D_T$ are the longitudinal and transverse dispersion coefficients.
We do not use different notations for the one- and two dimensional
particle densities. We assume that their distinction will be clear
from the context.

It appears that the effect of transverse dispersion is much more
dramatic than that of longitudinal dispersion, while its solution
is much easier to derive. Therefore, we analyze longitudinal and
transverse dispersion separately and start by including transverse
dispersion first.

\subsection{Transverse dispersion} We use an approach proposed
by~\cite{VanKoot}. Consider two distinct species, one adsorbing and one
non-adsorbing. Let the spatial distribution of the non-adsorbing
solute be given by $c(x,y,t)$. Further, let $\tau$ be the `free
residence time' of the adsorbing particles and let the
distribution of $\tau$ at time $t$ be $h_{ij}(\tau, t)$, i.e. for
particles in phase $j $ with an initial unit pulse in phase $i$
(see by (\ref {Gidtau})). If the initial pulse is $N^{0}_i$, the
 particle fraction with free residence time $\tau$ at time $t$
becomes $ N^{0}_i h_{ij}(\tau, t )$, where $\tau < t$. The spatial
distribution of this fraction is equal to that of the
non-adsorbing particles at $t = \tau$, or:
$
\md N_{ij}(x,y,t) = N^{0}_i h_{ij}(\tau,t)c(x,y,\tau)\md\tau.
$

\noindent Summing fractions with $\tau$, $ 0 \le \tau  \le t$, we obtain:
\begin{equation}
N_{ij}(x,y,t) = N^{0}_i \int_0^t
h_{ij}(\tau,t)c(x,y,\tau)\,\md\tau. \label{methodeKoot}
\end{equation}
We apply Van Kooten's approach first to the case with only
transverse dispersion. For a non-adsorbing solute with advection
in the $x$-direction and dispersion in the $y$-direction ($D_L$ is
assumed zero) the distribution is:
\begin{equation*}
c(x,y,t) = \frac{1}{{2\sqrt {\pi D_T t} }}\exp\! \left( \!{ -
\frac{{y^2 }}{{4D_T t}}} \right)\delta(x-vt).
\end{equation*}
Insert this function in the integral (\ref{methodeKoot}). Because
of the delta function the integral can be evaluated directly.
After substitution of $\tau$ by $x/v$ we obtain:
\begin{equation}
N_{ij} (x,y,t) = N^{0}_i
    \frac{ h_{ij} (\frac{x}{v} ,t)}{2{\sqrt {\pi D_T \frac{x}{v}} }}
    \exp\!\! \left( { - \frac{{y^2 }}{{4D_T \frac{x}{v}}}} \right).
\label{2dsolution}
\end{equation}
For a given value of $x$, expression (\ref{2dsolution}) describes
the distribution in $y$-direction of a certain amount of
particles. The amount is equal to $N^{0}_i h_{ij} (\frac{x}{v}
,t)$ per unit of length in $x$-direction and it spreads in the
$y$-direction as a Gaussian distribution with variance $2D_T
\frac{x}{v}$. The transverse variance now depends on the
$x$-coordinate, which is clearly in conflict with a 2D Gaussian
distribution. The interdependence of transverse variance and
$x$-coordinate can be understood by considering the residence
    times in the free phase. During the `free phase' time particles travel in
    positive $x$-direction and simultaneously
    spread in the $y$-direction. At a given time $t$ the particles that have spent more time in the free phase are found further along
    the $x$-direction. They are also more widely spread in
    $y$-direction, since they have been subject to dispersion for a longer time. Further, note that for $x\over v$$>t$ the functions $h_{ij} (\frac{x}{v} ,t)$
are zero.  Therefore, $N_{ij} (x,y,t) $ is zero for $x >vt$.
\begin{figure}[t]
\centering
\includegraphics[width =10cm]{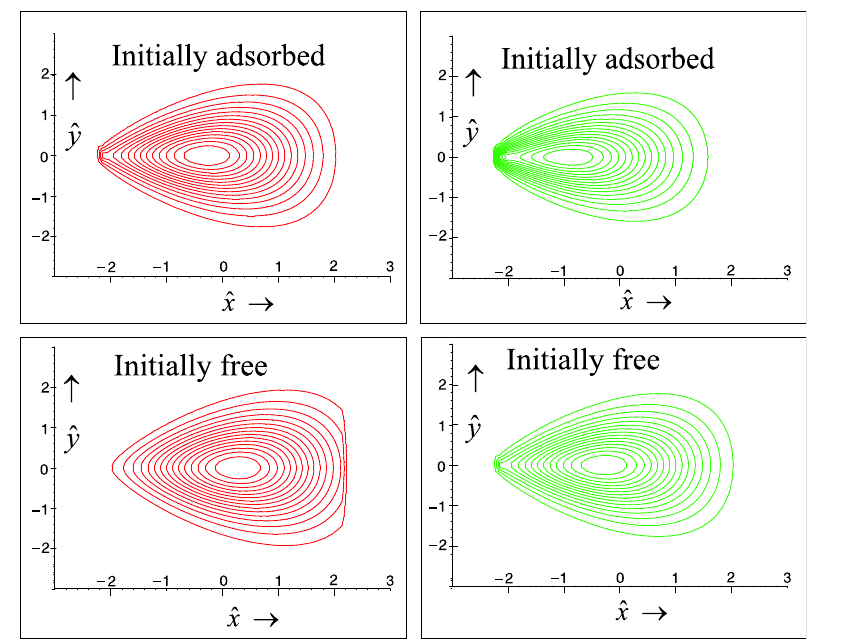}
\caption{Contour-lines for free phase (red) and adsorbed phase
(green) for several initial conditions. The symmetric case
($\lambda = 5$ and $\mu = 5$) at time $t$
=1.}\label{fig:onlytransverse}
\end{figure}

In Figure \ref{fig:onlytransverse} this is illustrated by 2D
contours for unit initial pulses in the free
 and  adsorbed phase. The following scaled
coordinates are used:
\begin {equation}
\hat x = \frac{{x - v^{*} t}}{{\sqrt {2D^{^{*}} t} }};  \quad\hat
y = \frac{y}{\sqrt {2D_T t} }.\label{scales}
\end {equation}
For a Gaussian distribution with dispersion coefficients $D^*$ and
$D_T/R$ in the longitudinal and transverse direction one would
expect elliptic contour-lines. Figure \ref{fig:onlytransverse},
however, shows typical `cigar'-shaped contours. It is clear that
close to the original injection point (in Figure
\ref{fig:onlytransverse} around $\hat x = -2$) the particles have
spend a relatively short time, $\tau$, in the `free' phase. During
this time the particles displace not only very little in
longitudinal direction, but also the spread in transverse
direction is very limited. The result of this is that $x$
displacement and $y$-spread are no longer independent, which is in
contradiction to a truly Gaussian system. Also note that in the 1D
case in Figure \ref{figureno1} the same parameters have been used
and the corresponding distribution in $x$-direction comes out
almost as a Gaussian curve (Figure \ref{figureno1}, upper right
graph).
\subsection{Longitudinal and transverse dispersion}

Analytical solutions including longitudinal and transverse
dispersion have been obtained by \cite{Carnahan}, but not in a
closed-form. The integral expression by Carnahan is, in essence,
equal to the one obtained by Van Kooten's approach:
\begin{equation}
N_{ij}(x,y,t) = \!\frac{N^{0}_i}{\sqrt {D_L D_T }}\!\int_0^t
\!{\frac{h_{ij} (\tau,t)}{{4\pi \tau }} \me^{{ - \frac{{(x - v\tau
)^2 }}{{4D_L \tau }} - \frac{{y^2 }}{{4D_T \tau }}}}
\md\tau}\label{koot2D}
\end{equation}
We apply expression (\ref{koot2D}) to evaluate the `full' system
(i.e, including advection, dispersion and sorption) using the
following values: $\lambda$ = 0.2, $\mu$ = 0, $v$ = 1, $D_L$ = 0.5
and $D_T$ = 0.1. The kinetics induced dispersion coefficient $D^*$
now becomes 0.62, which is in the same order of magnitude as
$D_L$.

\begin{figure*}[t]
\centering
\includegraphics[width =12cm]{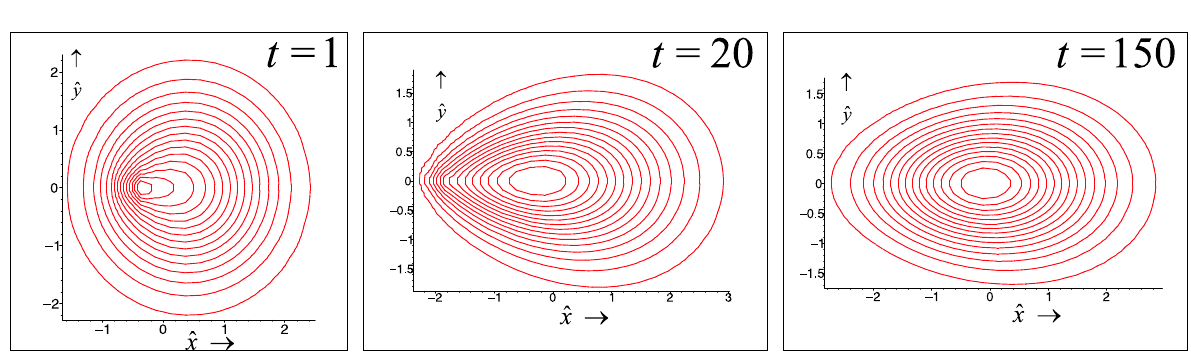}
\caption{Full two-dimensional case. Distribution of the free phase
at various times for an initial pulse in equilibrium; symmetric
$\lambda = 0.2$ and $\mu = 0.2$; $D_L$ = 0.5; $D_T$ = 0.1; $v$ = 1
} \label{fig:full2D}
\end{figure*}

The contour-lines in Figure \ref{fig:full2D} represent the
distribution of the free phase for $t$ = 1, $t$ = 20 and $t$ =
150, for the case of an initial pulse in equilibrium. Note that
for longitudinal and transverse coordinate we have applied the
same scaling as in Figure 4, (see (\ref{scales})). At a short time
($t$ = 1) almost circular contours occur, which suggests
Gaussianity. At this stage microdispersion dominates the spreading
process, which now progresses in a Gaussian way. At an
intermediate time ($t$ = 20) the `cigar'-shaped contours start to
develop. Here, we observe an increasing influence of the non
Gaussian kinetics-induced dispersion. Finally, at large times ($t$
= 150) the distribution becomes Gaussian again, with elliptic
contours. The parameters of the early and late Gaussian
distribution are quite different. For short times we have velocity
$v$ and dispersion coefficients $D_L$ and $D_T$. For large times
the velocity becomes $v^*$ and dispersion coefficients $ D_L/R +
D^*$ and  $ D_T/R$ for the longitudinal and transverse direction
respectively. The early Gaussian distribution occurs because
directly after the start of the pulse the effect of kinetic
exchange between the phases is still small and the free phase
consists mainly of particles that did not yet change their state.
Therefore, they behave as a non-adsorbing solute. As time goes on
the influence of kinetics becomes more apparent and `cigar' shaped
contours develop.
\section{Moment Analysis}

\subsection{One-dimensional}\label{1Dmoment}

The method of moments was introduced by \cite{Aris}, while
studying the flow of solutes through a capillary tube (Taylor
dispersion). Later, the method was successfully applied to solute
transport in a layered porous medium by \cite{Marle}.
 For contaminant transport with sorption the method was applied by
\cite{ValJH}, \cite{Valocchi}, \cite{Chrysikopoulos} and
\cite{Souadnia}. These authors examine the asymptotic values of
zero, first and second moments. Solutions for the entire
pre-asymptotic regime of zeroth, first and second moments are
derived by \cite{Michalak}. However, one of their expressions for
the central second moment turns out to be incorrect. We copy here
the formula (the $\oplus$ is a $+$), but should be a $-$) from
\cite{Michalak}, page 2136 for the normalized central moment,
where the solute is in the free phase at time 0 \emph{and} at time
$t$ :
\begin{equation*}\label{eq:Michalak}
\begin{split}
\sigma_{ff}^2&=\!{\frac {{t}^{2}A{v}^{2}\beta\, \left( \beta-1
\right) ^{2}}{ \left( \beta+1 \right) ^{2} \left( 1+\beta\,A
\right) ^{2}}}+ \! t \!\left( { \frac {2{ D}}{\beta+1}}+{\frac
{2{v}^{2}\beta}{k \left( \beta+1
 \right) ^{3}\!}} \right)\\
 &\quad + t A \left( {\frac {4{v}^{2}\beta\, \left( -{
\beta}^{2}A-{\beta}^{2}-\beta+1 \right) }{k \left( 1+\beta\,A
\right) ^{2} \left( \beta+1 \right) ^{3}}}\right)  + t A \left({\frac {2{\it
D}\,\beta\,
 \left( \beta-1 \right) }{ \left( \beta+1 \right)  \left( 1+\beta\,A
 \right) }} \right) \\
 &\qquad + {\frac {2{v}^{2}\beta\, \left( 1-A \right)
 \left( 3\,{\beta}^{2}A-3-\beta\, ( A\oplus 1 )  \right) }{{k}^{2
} \left( 1+\beta\,A \right) ^{2} \left( \beta+1 \right) ^{4}}}
+{\frac {4{\it D}\,\beta\, \left( 1-A \right) }{k \left(
1+\beta\,A
 \right)  \left( \beta+1 \right) ^{2}}}.
 \end{split}
\end{equation*}
Here Michalak and Kitanidis abbreviate
$A=A(t)=\exp(-(\beta+1)kt)$, and use the notation
$$ k=\mu,\quad \beta=\lambda/\mu.$$
By giving the expressions for $ \sigma _{ff}^2 $, $ \sigma _{fa}^2
$, $ \sigma _{af}^2 $ and $ \sigma _{aa}^2 $ Michalak and
Kitanides suggest that these expressions may be used to obtain
variances for a general initial condition, by linearity. However,
such a superposition can only be composed for the non-centralized
moments. At large times the first and second moment appear to
increase at a constant rate, which suggest that there exist an
effective velocity $v_{e}$ and effective dispersion coefficient
$D_{e}$:
$$
v_{e} =\frac{\mu}{{\lambda+ \mu}} v,\qquad D_{e} =
\frac{{\lambda\mu}}{{(\lambda+ \mu)^{3} }}v^{2} +
\frac{\mu}{{\lambda+ \mu}}D.
$$
\subsection{Two-dimensional}
The 2D case is described by the equations (\ref{eq2D-2}). We may
distinguish two different type of moments. The first type is that
of the marginal moments, or moments with respect to the
$x$-coordinate ignoring the information on the $y$-coordinates of
the particles. The second type consist of conditional moments,
either a moment with respect to $x$ for a given value of $y$, or a
moment with respect to the $y$ and for a given value of $x$.
Interestingly, the marginal moments for the 2D case are identical
to the moments for the 1D case. In the following subsections we
discuss the conditional moments.

\subsubsection{The $x$-moments conditioned on $y$}

This category of moments represent expected values of $x^{n}$ for
a population of particles with a specific $y$-coordinate. Note
that these moments are a function of $y$.  We define these moments
as:
\begin{equation*}
{M}_{f}^{(n)}\left(  y\right) =\int\limits_{-\infty}^{+\infty}{x^{n}N_{f}(x,y)\,\mathrm{d}x},
\,\, \text{ and }
{M}_{a}^{(n)}\left(  y\right) =\int\limits_{-\infty}^{+\infty}{x^{n}N_{a}(x,y)\,\mathrm{d}x}
\end{equation*}
For $n=0$ the expressions represent for each phase the particle
distribution along the $y$-direction (the total mass of particles
with the specified $y$-coordinate). Note that the higher order
moments ($n > 0$) are not yet divided by the zeroth moment, so
higher order moments are not normalized. For the initial condition
we consider a (unit) pulse with the phases in equilibrium. Then,
the two-dimensional particle distribution becomes:
\begin{equation}
N_{i}^{eq}(x,y) = \!\int_0^t \!\!h_{i}^{eq} (\tau,t){\frac{
 \me^{{ - \frac{{(x - v\tau )^2 }}{{4D_L \tau }} - \frac{{y^2
}}{{4D_T \tau }}}} }{{4\pi \tau \sqrt {D_L D_T }}}
\md\tau}\label{Kooteq}
\end{equation}
where $i$ is  $a$ or $f$. This expression is obtained by applying
(\ref{koot2D}), replace $h_{ij}$ by $h_{i}^{eq}$ (see
(\ref{eq_definition})) and take $N^0_i$ equal to 1. For the
equilibrium initial condition the moments are:
\begin{equation*}
\begin{split}
{M}_{i}^{(n)}\left(  y\right)    &
=\int\limits_{-\infty}^{+\infty}{x^{n}N^{eq}_{i}(x,y)\,\mathrm{d}x}
\end{split}
\end{equation*}
When we use (\ref{Kooteq}) and change the order of integration,
the zeroth, first
and second moments become:\\
\begin{equation}
M_{i}^{(0)}(y)  = \!\int_{0}^{t} h_{i}^{eq}
(\tau,t) \frac{\me^{  {-\frac{{y^{2}}}{{4D_{T}\tau}}%
}}}{2\sqrt{\pi D_{T}\,\tau}}  \mathrm{d}\tau\label{zeroMom}
\end{equation}
\begin{equation}
M_{i}^{(1)}(y)   = \!\int_{0}^{t}\!\! v \tau \, h_{i}^{eq}
(\tau,t) \frac{\me^{  {-\frac{{y^{2}}}{{4D_{T}\tau}}%
}}}{2\sqrt{\pi D_{T}\,\tau}}  \mathrm{d}\tau %
\label{firstMom}
\end{equation}
\begin{equation}
M_{i}^{(2)}(y) \!\!=  \!\!\!\int_{0}^{t}\!\!
(2D_{L}\tau\!+v^{2}\tau^{2})h_{i}^{eq}
(\tau,t) \frac{\me^{  {-\frac{{y^{2}}}{{4D_{T}\tau}}%
}}}{2\sqrt{\pi D_{T}\tau}}  \mathrm{d}\tau.
\end{equation}
\begin{figure*}
\centering
\includegraphics[width =12cm]{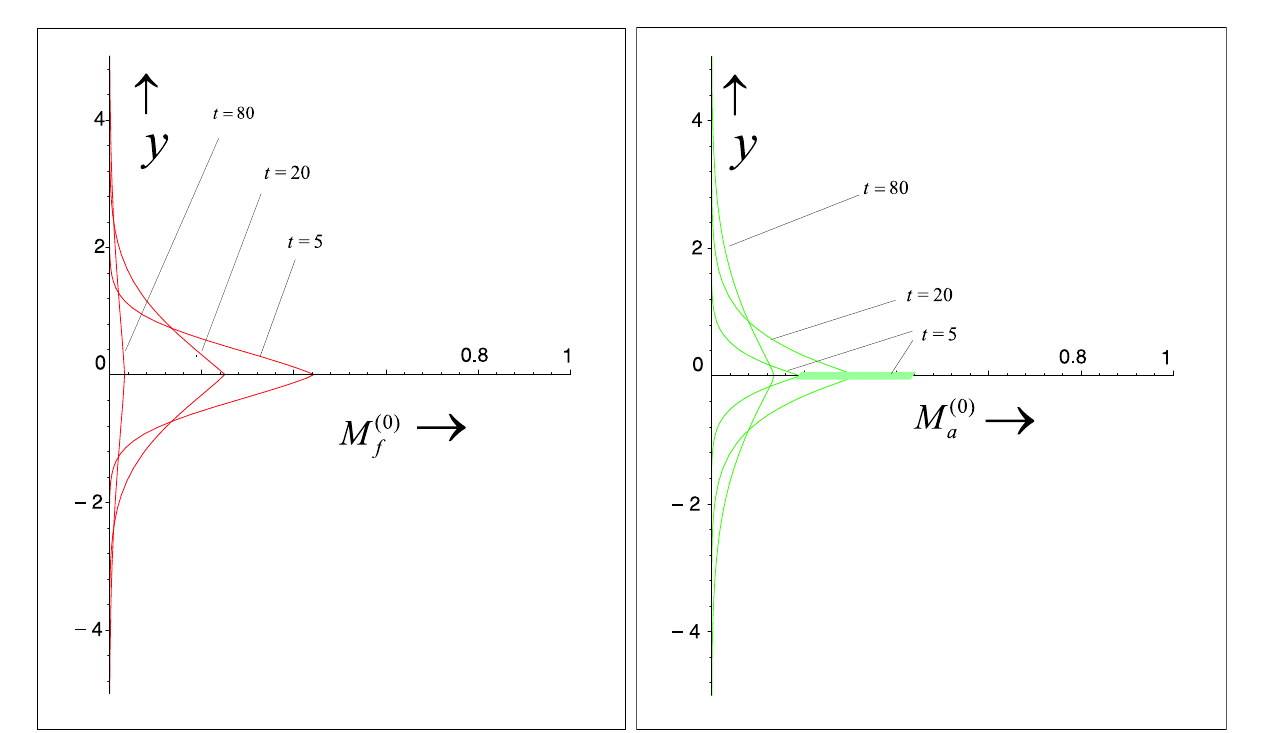}
\caption{Zeroth moments as given by
(\ref{zeroMom}) versus $y$ for $t$=5,20,80, left free phase, right
adsorbed phase; $\lambda$ = $\mu$ = 0.05; $D_{T}$ =0.03. Zeroth
moment on the horizontal axis. In the figure at the right hand
side the curve for $t$=5 has a pulse-shaped component $\delta(y)$
at $y$=0.}\label{fig:mom0a}
\end{figure*}
\\Figure \ref{fig:mom0a} shows the distribution along the $y$-axis
for the zeroth moments of free and adsorbed phase for several
times. The initial condition here is the equilibrium situation.
The zeroth moment is plotted horizontally. The figure shows that
particles of both phases gradually spread out in $y$-direction.
Note in the figure at the right hand side that at $t$=5 the
adsorbed phase is concentrated along $y=0$ and the curve has a
pulse-shaped component $\delta(y)$ for $y$=0. This pulse
represents the particles that did not yet spend time in the free
phase and are still in the initial position. At larger times most
adsorbed particles do have spent time in the free phase and the
spreading in $y$-direction becomes visible.

\begin{figure}[h] \centering
\includegraphics[width =8cm]{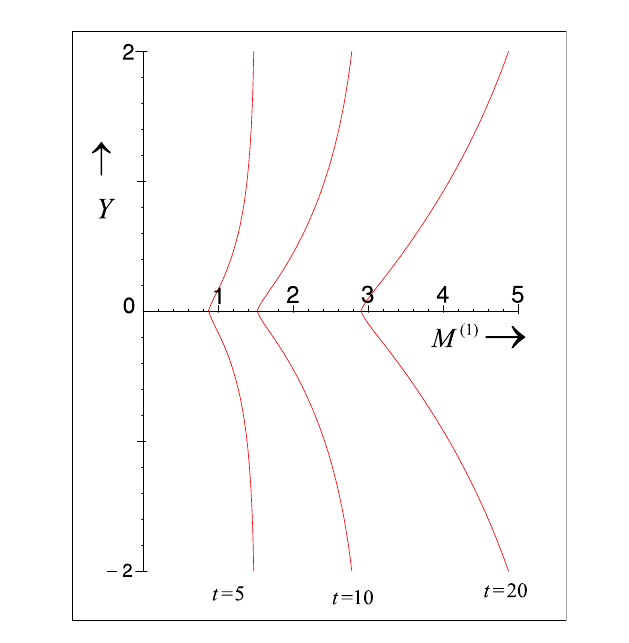}
\caption{First moments for free phase versus $y$
from (\ref{firstMom}) for $t$ $=5, 10, 20$. Numerical values:
$\lambda= 0.05$; $\mu=0.05$; $v = 0.3$; $D_{L}= 0.3$;
$D_{T}=0.03$. First moment on the horizontal axis.}%
\label{fig:mom1a}%
\end{figure}

Figure \ref{fig:mom1a} shows the first moment for the free phase.
The first moment, plotted horizontally, can be interpreted as the
average distance travelled in the $x$-direction for the particles
with a specific $y$-coordinate. As it appears, for $y$ = 0 the
first moments are smaller than for values of $y$ greater or
smaller than $0$. It is an alternative illustration of the tailing
effect. Particles spending less time in the free phase have less
opportunity to displace in the $x$-direction and spread in the
$y$-direction. Thus, a considerable fraction of these particles
are found around the point of origin. Further away from the
$y$-axis particles occur that have been able to disperse
laterally.  Therefore, they did spend some time in the free phase
and, consequently, were also displaced further along the $x$-axis.
Once more, we conclude that displacements in $x$- and
$y$-direction are mutually dependent. This dependence in the
spreading pattern is non-Gaussian.
\subsubsection{The $y$-moments, conditioned on $x$}

The $y$-moments conditioned on $x$ are defined as:
\begin{equation*}
{M}_{f}^{(m)}(x)={\displaystyle\int\limits_{-\infty }^{\infty}}
y^{m}{N}_{f}\left(  x,y\right)  \,\text{d}y ,\quad \text{and}\quad
{M}_{a}^{(m)}(x) ={\displaystyle\int\limits_{-\infty }^{\infty}}
y^{m}{N}_{a}\left(  x,y\right)  \, \text{d}y,
\label{x-moments}
\end{equation*}
We analyze the $y$-moments by the Laplace Transform using the
following initial conditions:
\begin{equation*}
\left. {N_f(x,y,t)}\right|_{t = 0} = \pi_f \delta(x)\delta(y),\quad \text{and}\quad
\left. {N_a (x,y,t)}\right|_{t = 0} = \pi _a \delta (x)\delta
(y)
\end{equation*}
After applying the Laplace transform to eq (\ref{eq2D-2}) and
eliminating the transform of $N_a$ we obtain the following
differential equation:
\begin{equation*}
D_L \frac{{\partial ^2 \widehat N_f }}{{\partial x^2 }} + D_T
\frac{{\partial ^2 \widehat N_f }}{{\partial y^2 }} -
v\frac{{\partial \widehat N_f }}{{\partial x}} - sb\widehat N_f  +
b \pi_f \delta (x)\delta (y) = 0
\end{equation*}
where $\widehat{N}_f$ denotes the particle density in Laplace
space, $s$ is the Laplace parameter and $b$ is an expression
depending on $s$:
\[
b = \frac{s + \lambda +\mu}{s + \mu}
\]
$\widehat N_a$ is related to $\widehat N_f$ by (see eq.~(\ref{eq2D-2})):
$
(s + \mu) \widehat N_a  - \pi _a \delta (x)\delta (y) = \lambda
\widehat N_f .
$

\begin{figure}[h]
\centering
\includegraphics[width =9cm]{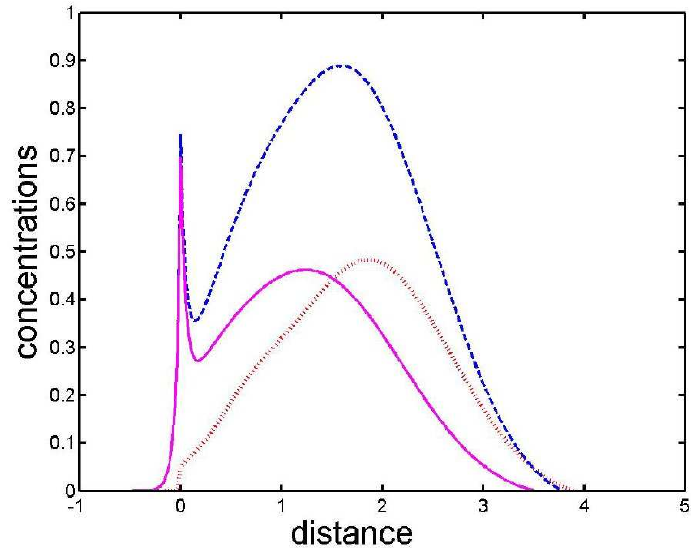}
\caption{$Z$eroth moment for
total concentration $M^{(0)}  $ (dashed), free concentration
$M_{f}^{(0)}$ (dotted)\ and adsorbed $M_{a}^{(0)} $ concentration
(drawn). Small numerical inaccuracies lead to the fact that the
zeroth moment for the total concentration is below the curve for
the free concentration. The peaks at $x=0$, representing the
function $\delta\left(  x\right) $ are \textquotedblleft
approximated" by
two exponential curves for reasons of illustration.}%
\label{mzeroy}%
\end{figure}

After taking the moments with respect to $y$ we find for the
zeroth moment:
\begin{equation*}
\begin{split}
& D_L \frac{{\partial ^2 \widehat M_f^{(0)} }}{{\partial x^2 }} -
v\frac{{\partial \widehat M_f^{(0)} }}{{\partial x}} - sb\widehat
M_f^{(0)} + b \pi_f \delta (x) = 0\\
& \widehat{M}_{a}^{(0)}
=\frac{\lambda}{s+\mu}\widehat{M}_{f}^{(0)}+ \frac
{\pi_a}{s+\mu}\delta(x)
\end{split}
\end{equation*}
The solution for $\widehat{M}_{f}^{(0)}=\widehat{M}_{f}^{(0)}(x)$ is 
\begin{equation*}
\begin{split}
&\widehat{M}_{f}^{(0)} =  \frac{{b\pi _f }}{{\sqrt {v^2 + 4bsD_L }
}}\exp\!\left(\!\frac {xv-\left|x\right|\!\sqrt {v^2 + 4bsD_L  } }{2D_L} \right)\\
\end{split}
\end{equation*}
where we use that the concentrations are zero at $x=\pm\infty.$
The inverse can be found with numerical Laplace inversion
algorithms \cite{davies}, \cite{stehfest}, \cite{valko}, and is
shown in Figure \ref{mzeroy}.

For the second moment we obtain the equations:
\begin{equation*}
D_L \frac{{\partial ^2 \widehat M_f^{(2)} }}{{\partial x^2 }} -
v\frac{{\partial \widehat M_f^{(2)} }}{{\partial x}} - sb\widehat
M_f^{(2)} + 2
D_T \widehat M_f^{(0)} = 0, \,\, \text{and}\,\,
s\widehat{M}_{a}^{(2)}   =\lambda \widehat{M}_{f}^{(2)}-\mu
\widehat{M}_{a} ^{(2)}
\end{equation*}
\noindent where we use that
$
{\displaystyle\int\limits_{-\infty}^{\infty}}y^{2}\frac{\partial^{2}\widehat{N}_{f}%
}{\partial y^{2}}\,\mathrm{d}y=\widehat{M}_{f}^{(0)}
$,
with the solution:%
$$
\widehat{M}_{f}^{(2)}(x) = b \pi_f \frac{\mu^2 D_L}{\lambda +
\mu}\frac{\exp \left( \frac {xv -  x \sqrt {v^2 + 4bsD_L  }
}{2D_L}
\right) } {v^2 +4bsD_L} \times \frac{\mu}{v}\left(
x+\frac{2 D_L}{\sqrt{v^2 +4bs D_L } }\right),\,\, \text{ for $x>0$},
$$
and $\widehat{M}_{a}^{(2)}(x) ={\displaystyle \frac{\lambda}{s+\mu}\widehat{M}_{f}^{(2)}(x)}$.

The ratio of the numerically~\cite{stehfest} inverted $\widehat{M}_{f}%
^{(2)}$ and  $\widehat{M}_{f}^{(0)} $ is the normalized and
centralized second moment (conditioned on $x$) and this ratio
turns out to be $x/P\!e$. Thus, the lateral spread depends on $x
$, which (again) explains the tailing and the previously observed
`cigar shape'.

\section{Conclusions}
We discuss a stochastic particle approach, considering a particle
that changes between a mobile and immobile state governed by a
Markov chain, while its spatial displacement is governed by a
random walk. Generally, particle models are assumed to describe
the advective-dispersive-kinetic transport process correctly, but
in the literature a rigorous proof is lacking. Our analysis shows
that it is possible to derive the correct set of differential
equations from a stochastic model for a single particle.\\
To examine the non-Gaussian nature of the spreading process we
 analyze first the telegrapher's equation. This equation arises
when only advection and kinetic sorption is considered. We show
that in such a system an apparent dispersion process occurs,
generated only by the kinetic changing of the particle in states
with different velocities (i.e. $zero$ or $v$). This
`kinetics-induced' dispersion is non-Gaussian for short and
intermediate times, while at large times the process develops as
Gaussian dispersion. When hydrodynamic dispersion is included
again, the spreading process becomes a combination of a Gaussian
and non-Gaussian dispersion. We illustrate this for the 2D case.
At short times the process is Gaussian, since hydrodynamic
dispersion is the dominating process. At intermediate times the
influence of kinetic induced dispersion increases and the
spreading becomes non-Gaussian. Finally, at large times the
dispersion becomes Gaussian again, but the (effective)
longitudinal dispersion coefficient has an additional term due to
the kinetics. Moveover, we find for the 2D case that the
transverse spreading depends on the longitudinal coordinate,
resulting in `cigar-shaped' contours. The mechanism is best
illustrated when longitudinal dispersion is assumed zero. Here the
particles displace in the $x$-direction by advection and spread
transversely by dispersion. Particles spending more time in the
adsorbed phase are displaced less in $x$-direction, and also less
spread out in $y$-direction. In a truly Gaussian distribution the
transverse spreading is independent of the longitudinal
coordinate. When longitudinal dispersion is included the same
effect is observed, although for short times (compared to the
kinetic exchange rate) the situation is now dominated by
hydrodynamic dispersion. With respect of the validity of effective
properties (velocity and dispersion), we conclude the following.
The velocity and dispersion coefficients are represented by the
rate of increase of the first and centralized second moment (times
$1/2$). For cases with low adsorption and desorption rates, the
rates of increase for the moments remain time dependent for a
relatively long time. We conclude that constant effective
properties can not be defined directly after the start of solute
injection. For large times an asymptotic behavior is observed with
a constant mean displacement and rate of dispersion, while the
third moment vanishes and the kurtosis approaches a value of 3.
This can be proved via our stochastic model by applying a
sophisticated version of the Central Limit Theorem (Section
\ref{MandV}). The critical time, required before dispersion
coefficients becomes constant, is in the order of $3/( \lambda+\mu
)$. The effective velocity is $v\mu/(\lambda+\mu)$, similar to the
case of linear equilibrium adsorption. However, at early times the
free particles move with the original groundwater velocity and
build up a lead with respect to the adsorbed phase. In the
asymptotic stage, the adsorbed and free particles displace with
the same average velocity, but the lead of the free plume is
maintained. The effective longitudinal dispersion, can be much
higher than in the case of linear equilibrium adsorption. The
additional term, $v^2 \lambda \mu /(\lambda  + \mu )^3$ depends
also on the groundwater velocity. It is remarkable that it depends
on $v^2$, while the micro-scale dispersion is linear in $v$.

\end{document}